\documentclass[conference,10pt]{IEEEtran_ICRATmod}

\IEEEoverridecommandlockouts
\usepackage[numbers]{natbib}

\usepackage{graphics} 
\usepackage{epsfig} 
\usepackage{times} 
\usepackage{amsmath} 
\usepackage{amssymb}  
\usepackage{amsthm}
\usepackage{amsfonts}
\usepackage{mathtools}
\usepackage{multirow} 
\usepackage{rotating} 
\usepackage{graphicx}
\usepackage{subfig}
\usepackage{float}
\usepackage{flushend}
\usepackage[table,xcdraw]{xcolor}
\usepackage{csquotes}
\usepackage{bm}
\usepackage{relsize}
\usepackage{epspdfconversion}
\usepackage{colortbl}
\usepackage{lipsum}
\usepackage[mathscr]{eucal}
\usepackage{fancyhdr}
\usepackage{gensymb,bbm}

\usepackage{tikz}

\usepackage{algorithm}
\usepackage[algo2e]{algorithm2e}

\usepackage{bbm}

\usepackage{hyperref}

\newcommand{\norm}[1]{\left\lVert#1\right\rVert}

\newcommand{\abs}[1]{\left\lvert#1\right\rvert}

\DeclareMathAlphabet\mathbfcal{OMS}{cmsy}{b}{n}

\def\compactify{\itemsep=0pt \topsep=0pt \partopsep=0pt \parsep=0pt}
\newcommand{\compress}{\itemsep=0pt \topsep=0pt \partopsep=0pt \parsep=0pt \leftmargin=30pt \labelwidth=30pt}

\let\latexusecounter=\usecounter

\begin{document}
\title{Spaceport Facility Location Planning\\within the US National Airspace System}

\author{\scriptsize
	   	   \IEEEauthorblockN{Haochen Wu, Kevin R. Sun, Jackson A. Miller, Oliver Jia-Richards, Max Z. Li}
	   	   \IEEEauthorblockA{University of Michigan\\Ann Arbor, MI, USA\\ \{\href{mailto:haocwu@umich.edu}{haocwu}, \href{mailto:krsun@umich.edu}{krsun}, \href{mailto:jacksomi@umich.edu}{jacksomi}, \href{mailto:oliverjr@umich.edu}{oliverjr}, \href{mailto:maxzli@umich.edu}{maxzli}\}@umich.edu}
}

\vspace{-2cm}
\maketitle

\renewcommand{\headrulewidth}{0pt}
\lhead{\textcolor{black}{ICRAT 2024}}
\rhead{\textcolor{black}{Nanyang Technological University, Singapore}}
\cfoot{\textcolor{black}{\thepage}}
\lfoot{}
\thispagestyle{fancy}
\pagestyle{fancy}

\fancypagestyle{firststyle}
{
    \fancyhf{}
    \lhead{\textcolor{black}{ICRAT 2024}}
\rhead{\textcolor{black}{Nanyang Technological University, Singapore}}
    \cfoot{\textcolor{black}{\thepage}}
    \lfoot{\textcolor{white}{\thepage}\\ \vspace{0.1cm} \scriptsize \parbox{0.465\textwidth}{\hrule ~~\\ The authors gratefully acknowledge partial funding support from the University of Michigan Seeding To Accelerate Research Themes (START) program.}}
}

\noindent
\begin{abstract}
The burgeoning commercial space transportation industry necessitates an expansion of launch infrastructure to meet rising demands. However, future operations from these large-scale infrastructures can result in new impacts, particularly to air traffic operations. To rigorously reason about where such future spaceports might be located and what their impacts might be, we introduce a facility location planning model for future US spaceports (SPFLP). Central considerations for the SPFLP include population density, space launch trajectories, and potential impacts to air traffic within the US National Airspace System (NAS). The SPFLP outputs a cost-optimal set of candidate locations for future spaceports while satisfying a range of operational constraints. By conducting sensitivity analyses on the SPFLP, we are able to examine differences in flight rerouting costs and optimal launch mission allocations. Our model and numerical experiments offer valuable insights for future spaceport site selection, contributing to the strategic development of commercial space transportation while keeping in mind the need to integrate these operations within the NAS.
\end{abstract}
\begin{small}{{\bfseries\itshape Keywords---Air-space integration; Commercial space; Spaceports; Facility location planning}}\end{small}

\thispagestyle{firststyle}
\section{Introduction}  \label{sec:intro}
Space launches have become critical in numerous sectors, including military, navigation, weather forecasting, disaster response, and entertainment \cite{spaceport_development_2008}. Within this broad spectrum, commercial space launches have garnered significant attention: The Federal Aviation Administration's (FAA) Office of Commercial Space Transportation categorized several potential markets within the commercial suborbital space launch sector \cite{faa_suborbital}. 
Furthermore, with recent orbital launch successes demonstrated by SpaceX and other commercial launch vehicle companies, the barrier of entry into space continues to decrease while the demands for space exploration and commercial utilization continues to grow. The trend is evidenced by the steady increase in the number of FAA-licensed and permitted launch operations over the past seven years, with projections indicating an exponential rise through 2026 \cite{space_data,surging_launch}.

\subsection{Motivation and research problem}
As US commercial and military space interests mature and expand, its space infrastructure must keep pace to accommodate space launch demands. The establishment of new spaceports will eventually become vital to broadening access to space; current space launches mainly utilize federal ranges such as Cape Canaveral Space Force Station, Kennedy Space Center, and Vandenberg
Space Force Base. These federal range sites are often preferred due to their ability to service vertical launches 
and their more affordable price. 
In recent years, commercial spaceports like Spaceport America and Blue Origin's Launch Site One have also been built to offer more customized launch infrastructure in support of new missions~\cite{spaceport_development_2008}.

Our work herein is motivated by the challenge of determining \emph{where} future spaceports may be located, then examining the \emph{impact} launch operations from these future sites may have on the US National Airspace System (NAS). 
In particular, traditional traffic management approaches for managing interactions between  NAS volume and launch/reentry operations are becoming inefficient due to the rise in commercial space launches \cite{FAA_conops_integration}. The frequency of these launch/reentry operations is likely to cause substantial increases in flight times, delays, and fuel burn, culminating in escalated costs \cite{space_data}. Furthermore, these impacts may be inequitably felt throughout the NAS \cite{IntegratedSimulation}. Consequently, it will be critical to identify locations that not only meet the prerequisites for spaceports but also consider their potential impact on NAS flight operations.
To address this challenge, we utilize the framework offered by facility location planning problems \cite{daskin1997network}, tailoring it to specifically identify optimal locations for new spaceports in a systematic and data-driven way. 

\subsection{Background and prior works}
We first examine prior research on spaceport location selection. Ref. \cite{spaceport_development_2008} sets forth guidelines for both public and private commercial spaceports, covering aspects such as siting, layout, and supporting infrastructures. These guidelines suggest that spaceports be ideally situated in sparsely populated areas and away from flight corridors to allow for flexible launch schedules. Another study focuses on optimizing potential spaceport locations in Australia, identifying key factors such as population density, proximity to major cities, airspace considerations, thunderstorm activity, and cloud coverage \cite{wild2015optimising}. 
Furthermore, Refs. \cite{dachyar2018spaceport} and \cite{perwitasari2019indonesia} utilize the Analytic Hierarchy Process (AHP) for spaceport location determination. AHP is a multi-criteria decision-making process that emphasizes weighting criteria based on expert opinions. In Ref. \cite{dachyar2018spaceport}, factors like orbit type, transportation, flight trajectory, population density, and potential disaster risks are prioritized, while Ref. \cite{perwitasari2019indonesia}
focuses on flight trajectory, military sensitivity, area size, tracking, and mission objectives. These decision-making approaches are well-suited for long-term strategies involving a relatively smaller number of spaceports. However, given the anticipated substantial rise in space launch demand, there is a pressing need for more responsive and comprehensive spaceport location selection models. 

\subsection{Space launch operations within the NAS}

Refs. \cite{murray2008space,murray2009air} highlight the potential risks (e.g., debris from failed space launches) posed by space launch and reentry operations to air traffic, and propose methods (e.g., corridor-based safe zones for commercial flights) to mitigate these hazards. 
Ref. \cite{bilimoria2013space} provides detailed methodologies for generating space traffic corridors to separate spacecraft and flights, accounting for factors like target orbit inclination and launch azimuth. 
Additionally, Ref. \cite{IntegratedSimulation} develops an agent-based modeling and simulation framework to assess the impact of space launches on the air traffic system, revealing varying degrees of impact based on launch site location, launch slot, and failure probability. Finally, reinforcement learning-based methods have been used to study the problem of rerouting and metering aircraft during space launches \cite{tompa2018,tompa2019optimal}.

\subsection{Facility location planning (FLP) problems}
The FLP has been extensively studied and is applicable to network design and the location selection of critical infrastructures \cite{cornuejols1983uncapicitated}. FLP models aim to determine the locations of infrastructures which satisfy a set of constraints (e.g., demand) while minimizing various costs. 
The complexity of FLP models, especially due to the presence of integer variables and constraints, renders them to be NP-hard \cite{melkote2001capacitated}. To address this challenge, a linear programming relaxation technique with valid inequalities is proposed in Ref. \cite{melkote2001capacitated}, while Ref. \cite{tragantalerngsak2000exact} introduced an exact method utilizing a Lagrangian duality-based branch-and-bound algorithm. Spaceport location planning can be considered as a variant of the capacitated FLP model, where the capacities of the spaceports are limited, and the goal is to satisfy all space launch demands while minimizing cost.

\section{Contributions}
In this work, our contributions are as follows:
\begin{enumerate}
    \item We identify key factors influencing the selection of launch sites, develop an algorithm for computing feasible launch angles, classify launch mission types based on historical launch data, and estimate various types of associated costs including the transportation and operation cost, the launch cost, and the flight rerouting cost. 
    \item We introduce a spaceport location planning model that optimizes the selection of spaceport sites and mission assignment plans to achieve the lowest possible costs.
    \item We perform a sensitivity analysis to assess the impact of space launches on air traffic within the NAS and to identify ideal launch sites for different mission types.
\end{enumerate}
%

\section{Methodology}  \label{sec:methodology}
We first discuss the cost factors and constraints that go into the spaceport facility location planning problem (SPFLP), as well as construct a hypothetical but plausible portfolio of space launch missions that will be serviced by future spaceport locations (i.e., the \emph{demand} to be served by the \emph{supply} of selected spaceport locations). We then detail the formulation of the SPFLP. 

\subsection{Costs and demand}


We consider each US county as a potential candidate site for a spaceport and utilize county-level census data from the US Census Bureau \cite{us_census_2020} to estimate various parameters associated with hosting a spaceport, such as transportation cost and launch demand, as well as ascertain population densities to inform safety constraints. Transportation costs are defined as the expenses incurred to transport personnel and launch equipment to the spaceport. 
The operation costs of a spaceport are inclusive of both its construction expenses and the costs associated with spacecraft launch operations. Due to the lack of information for estimating such transportation and operation costs, in this paper we derive an approximation of the transportation and operation costs for each US county by analyzing the mean traffic commute times and the median house value as reported in the census data. 
Launch costs are attributed to the amount of fuel needed for launching a rocket with a particular mission profile. This depends on factors such as launch site geographic location, initial launch trajectory, and mission goals such as orbital inclination and semi-major axis; we detail this cost component in  
Section \ref{sssec:launc_traj_costs}. 

\begin{figure}[ht]
\centering
\includegraphics[width=.8\linewidth]{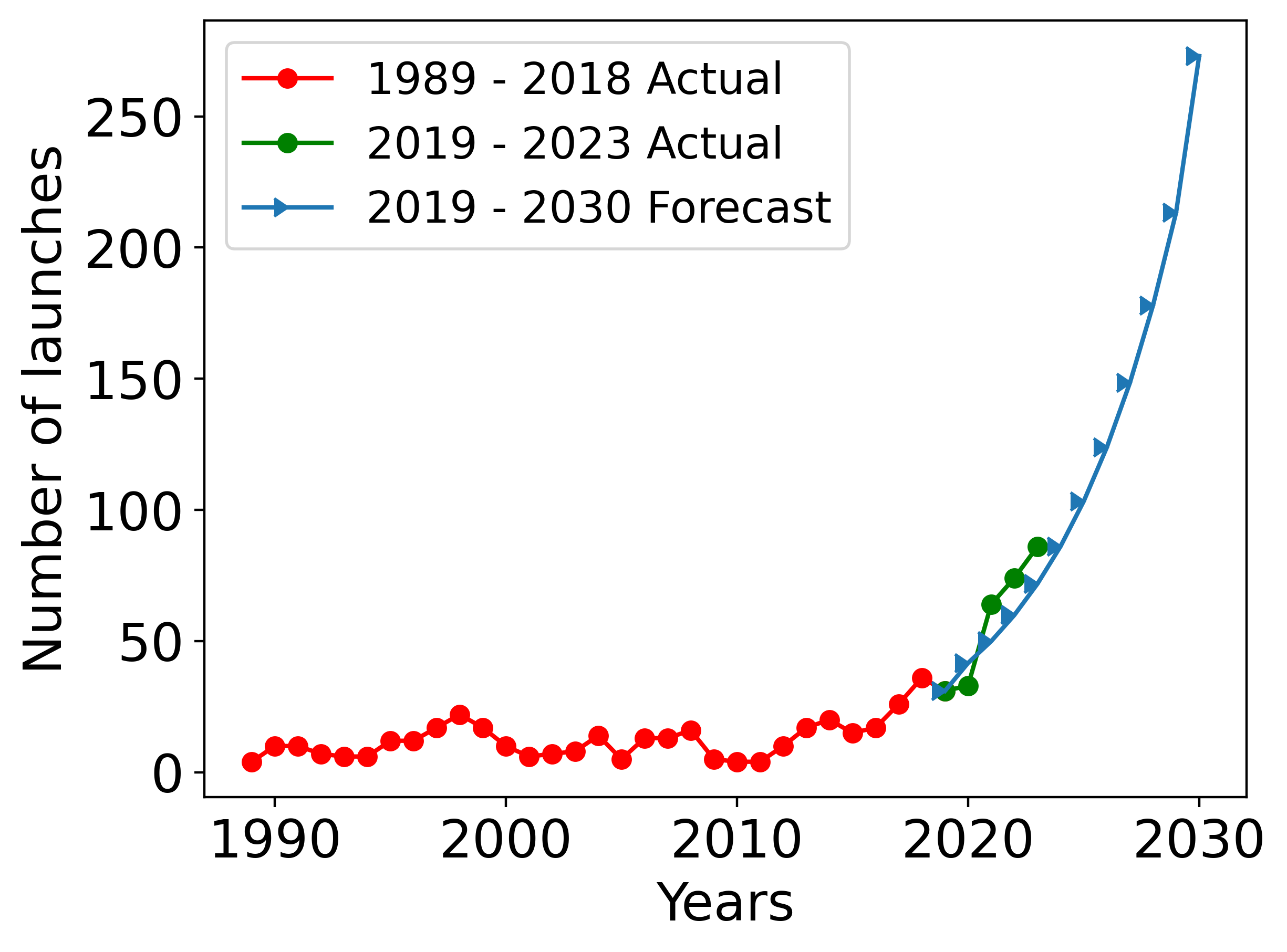}
\caption{Space launch demand forecast via \eqref{eq:forecast_eqns} with $\alpha = 0.6$ and $\beta = 0.05$.}
\label{fig:demand}
\end{figure}

A primary motivation for the SPFLP is to project potential future launch impacts to commercial aviation. 
Utilizing flight trajectory data from OpenSky \cite{schafer2014bringing}, we obtained the number of flight paths crossing through the airspace above potential spaceport locations. This information serves as a proxy through which \emph{flight rerouting costs} will be estimated. 
To estimate the total demand of space launch by 2030, we apply Holt's linear trend method from Ref. \cite{hyndman2018forecasting} on time series data of historical launches. This method, which is an extension of exponential smoothing, requires a forecast equation and two smoothing equations:
\begin{equation}
\begin{aligned}
\widehat{y}_{t+h} &= l_{t} + hb_{t},\\
l_t &= \alpha y_t + (1 - \alpha) \left(l_{t-1} + b_{t-1} \right),\\
b_t &= \beta\left(l_t - l_{t-1}\right) + (1-\beta)b_{t-1},
\label{eq:forecast_eqns}
\end{aligned}
\end{equation}

\noindent
where $t+h$ is the time period of interest, $l_t$ denotes an estimate of the level of the series (the average value in the series) at time $t$, $b_t$ denotes an estimate of the trend of the series (the average increasing or decreasing value in the series) at time $t$, and $\alpha \in [0,1]$, $\beta \in [0,1]$ are smoothing parameters for the level and trend, respectively. Using historical launch data from 1989 to 2023, we forecast via  \eqref{eq:forecast_eqns} the launch demand by 2030, with $\alpha = 0.6$ and $\beta = 0.05$. The result is shown in Figure \ref{fig:demand}. 
The estimated annual space launch demand for 2030 is 273, which will be taken as an input for the spaceport location planning model. Note that a more sophisticated scenario-based analysis could be conducted (e.g., low-, medium-, and high-demand forecast scenario for 2030). However, our focus in this paper is formulating and analyzing the results of the SPFLP.

%

\subsection{Launch trajectory constraints} \label{ssec:launch_trajectory}
To integrate launch trajectory constraints into the SPFLP, we first need to filter for all counties that can feasibly serve as the origin point of a rocket trajectory, obeying population-based safety constraints. Then, given the set of all feasible counties and launch angles, we detail the computation of the launch cost, which is specific to each county, desired mission type, and the launch angle. 




\subsubsection{Filtering feasible trajectories} \label{sssec:feasible_traj}

We first assume that the latitude and longitude coordinates for the centroid of a county, as given by the US Census Bureau data set \cite{us_census_2020}, will serve as the exact location of the future spaceport. From this latitude and longitude, we check the full range of launch azimuths $\psi \in [0\degree, 360\degree)$. The launch azimuth is the initial launch angle, measured clockwise from $0\degree$ North \cite{azimuth}. Note that we use \emph{launch angle} and \emph{initial launch azimuth} interchangeably.  For each angle, we form a hazard zone representing the area where populations and property could be harmed by debris from a launch failure.  If the population under the hazard zone was sufficiently low, we mark the tested azimuth as a feasible launch angle for a given county. Note that although the FAA dictates that launch operations must provide impact hazard areas, aircraft hazard areas, and ship hazard areas \cite{FAA_hazards}, for tractability and generalizability, we used generic hazard zones. Future work could refine and tailor these hazard zones to specific launch mission profiles and other needs. We also use these generic hazard zones to determine aircraft rerouting requirements.

For a given launch azimuth, we simplified the launch trajectory to be a line starting at a county's central latitude and longitude, extending outwards at angle $\psi$. To form the hazard zone, we then fan out boundary lines offset from the trajectory line at $\psi\pm\xi$ degrees, where $\xi$ is a \emph{buffer angle}. A larger $\xi$ indicates a larger buffer, i.e., a larger hazard zone.


The extent of these hazard zones depends on the time at which a rocket reaches orbit, since debris at that point are not considered harmful. As per the Orbital Degree Mitigation Standard Practices, the risk of human casualty from orbital debris impact should be less than 1 in 10,000 \cite{orbital_debris}. 
We conservatively truncate the hazard zone at the outer boundary lines of  the 2010 US Census Grid population map data set  \cite{population_dataset}. These restrictions effectively limit feasible launch trajectories to originate from coastal US counties, launched at an angle towards open waters. 
These restrictions conform to how current launches operate at sites such as Cape Canaveral and Vandenberg. Once a hazard zone is determined, we check for launch feasibility as determined by the population size within the hazard zone. We gathered feasible trajectories for counties under three possible hazard zone sizes: $\xi = 5\degree \text{(smallest)}, 7.5\degree,$ and $10\degree$ (largest). These feasible trajectory sets were used as three different inputs to the SPFLP model to test the effects of hazard zone size on model outputs. 

\subsubsection{Launch trajectory costs} \label{sssec:launc_traj_costs}
In order to quantify the insertion cost of launch mission type $j$ (which determines the target circular orbit with radius $r_j$ and the desired mission orbital inclination $\mathbbm{i}_j$) from a candidate spaceport location at county $i$ (with launch site latitude $\phi_i$ and feasible launch azimuth $\psi$), we calculate the propulsive cost of placing the rocket into the targeted orbit, and the propulsive cost to correct the orbital inclination. These costs are estimated as $\Delta v$, a standard measure of maneuver cost in astrodynamics.  
The total cost of inserting a rocket into the target circular orbit is $\Delta v = \Delta v_1 + \Delta v_2$, where $\Delta v_1$ is the launch cost and $\Delta v_2$ is the orbital inclination correction cost. We first examine the launch phase: the velocity of the rocket once it achieves orbit can be approximated as
\begin{equation}
\begin{aligned}
\Vec{v}^{\,} = 
\begin{bmatrix}
\Delta v^{2}_1 + 2\Delta v_1 \omega R \cos\phi_i \\
\Delta v_1 \cos\psi
\end{bmatrix},
\end{aligned}
\end{equation}

\noindent
where $\omega$ and $R$ are rotational velocity and radius of the Earth, respectively. The magnitude of the rocket's velocity is 
\begin{equation}
\begin{aligned}
\norm{\Vec{v}}^{2} = \Delta v^{2}_1 + 2 \Delta v_1 \omega R \sin (\psi) \cos (\phi_i) + \omega^{2}R^{2}\cos^{2}\phi_i.
\end{aligned}
\end{equation}

\noindent
This velocity magnitude must be equal to the desired orbit velocity. For a circular orbit around a planetary body with gravitational parameter of $\mu$, the required velocity for a circular orbit is $v_{o} = \sqrt{\mu / r_j}$. Therefore, we have that 
\begin{equation}
\begin{aligned}
\Delta v^{2}_1 + 2 \Delta v_1 \omega R \sin(\psi) \cos(\phi_i) + \omega^{2} R^{2} \cos^{2}\phi_i = \frac{\mu}{r_j},
\end{aligned}
\end{equation}
\noindent
and we can solve for $\Delta v_1$ as:
\begin{equation}
\begin{aligned}
\Delta v_1 = -\omega R \sin(\psi) \cos(\phi_i) + \sqrt{-\omega^{2} R^{2} \cos^{2}(\phi_i) + \frac{\mu}{r_j}},
\end{aligned}
\end{equation}
which is the launch cost neglecting gravitational and atmospheric drag losses as these are mostly independent of launch site latitude.

We can now calculate the initial orbital inclination $\mathbbm{i}_l$ that the rocket enters. Note that we may have that $\mathbbm{i}_l \neq \mathbbm{i}_j$, requiring a second maneuver to correct the orbital inclination. We first calculate the azimuth of the orbit velocity after launch
\begin{equation}
\begin{aligned}
\psi_{l} = \tan \left( \frac{\Delta v_{1} \sin(\psi) + \omega R \cos\phi_i}{\Delta v_{1} \cos\psi} \right),
\end{aligned}
\end{equation}
\noindent
then, the initial orbital inclination $\mathbbm{i}_l$ is given by 
\begin{equation}
\begin{aligned}
\mathbbm{i}_l = \arccos\left(\cos (\phi_i) \sin\psi_{l}\right). 
\end{aligned}
\end{equation}
\noindent
We can now calculate the orbital inclination correction cost $\Delta v_2$ to adjust $\mathbbm{i}_l$ to the desired mission inclination $\mathbbm{i}_j$ as
\begin{equation}
\begin{aligned}
\Delta v_2 = 2\sqrt{\frac{\mu}{r_j}} \sin\left(\frac{ \abs{ \mathbbm{i}_l- \mathbbm{i}_j } }{2}\right).
\end{aligned}
\end{equation}
\noindent
Note that $\abs{ \mathbbm{i}_l- \mathbbm{i}_j }$  represents the magnitude of the actual angular difference between the rocket's current inclination and the desired mission inclination. 

\subsection{SPFLP formulation}

We bring together the cost factors discussed previously in this section and formulate the full SPFLP model. Solving the SPFLP as a discrete optimization problem provides us not only with a set of candidate counties, but also an allocation of launch mission types for each county, which minimizes a total cost objective. We note that the feasible launch azimuths (i.e., feasible launch angles) for each county is pre-computed prior to the SPFLP (via Section \ref{sssec:feasible_traj}); the set of counties considered by the SPFLP are exclusively those with a non-empty set of feasible launch azimuths.





\subsubsection{Extracting and allocating mission types}
According to the FAA National Spaceport Network Development Plan \cite{FAA_spaceport_development_plan}, spaceports should be adaptable to either orbital or suborbital profiles as well as a variety of trajectories and inclinations. 
As discussed in Section \ref{sssec:launc_traj_costs}, the launch cost of a given mission from a given spaceport depends on the orbital characteristics of the mission as well as the geographic location of the spaceport. To address this within the SPFLP, we introduce an integer decision variable $y_{i,j}$ that accounts for the number of mission type $j$ assigned to a candidate spaceport at county $i$.  
The optimal allocation of missions to candidate spaceports can now be coupled with the decision of which counties to select for future spaceport locations. 

\begin{figure}[h]
\centering
\includegraphics[width=1\linewidth]{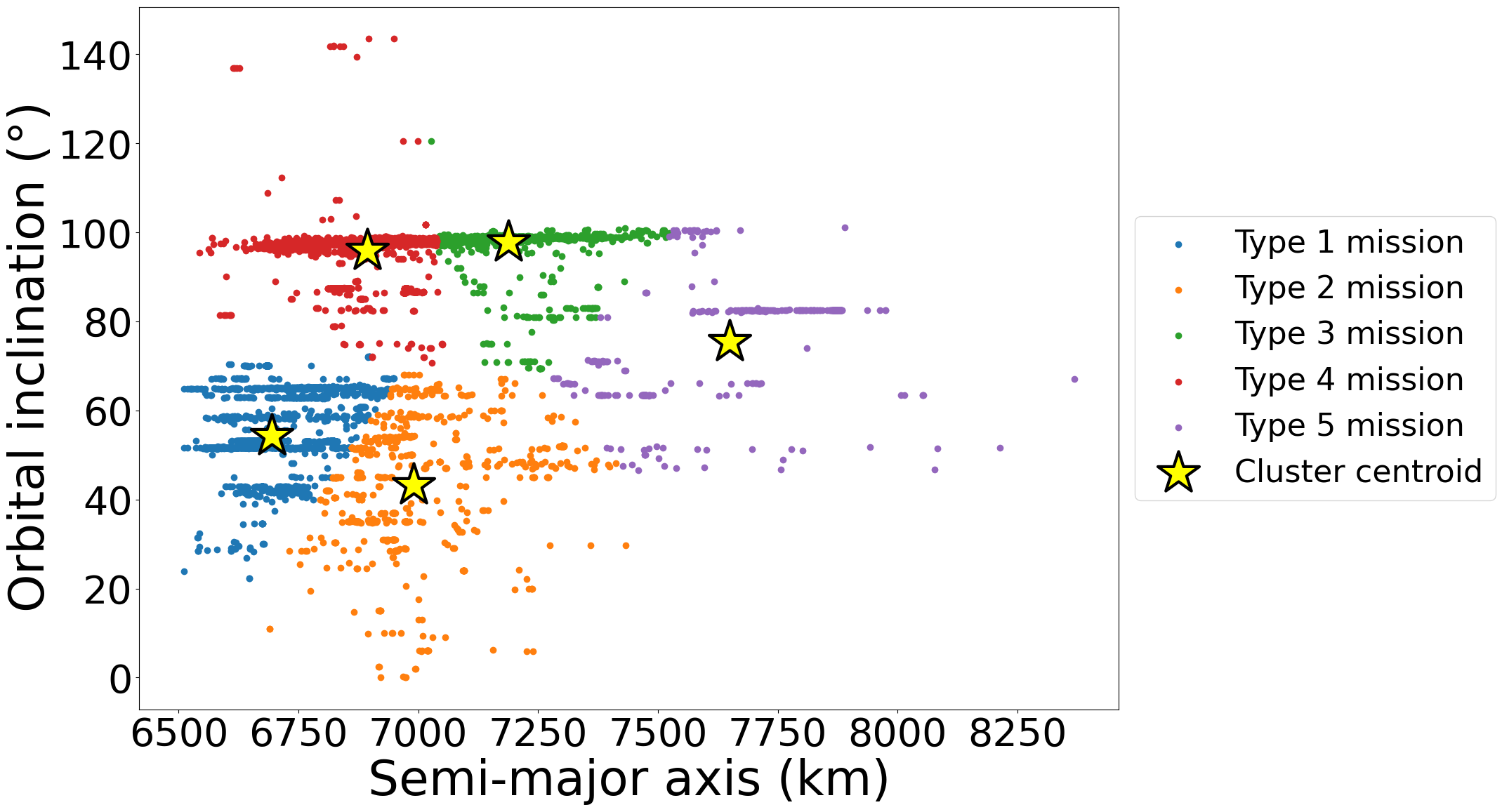}
\caption{Mission type clustering results based on historical orbital inclination and semi-major axis pairs.}
\label{fig:mission_cluster}
\end{figure}

Our approach for identifying different mission types is informed by historical launches. Specifically, we apply a simple clustering method ($k$-means) to a data set of historical launches \cite{orbit_data}. Note that spaceport stakeholders seeking to deploy our model should use their own launch mission portfolios---as this data is highly proprietary, we resort to historical data to extract mission types. From the data set of historical launches, each launch $\ell$ is represented by a data point $\left(a_\ell, \mathbbm{i}_\ell\right)$ comprising of the orbital semi-major axis $a_\ell$ and orbital inclination $\mathbbm{i}_\ell$. We cluster based on these two features using the Euclidean distance as the clustering metric. The Silhouette score \cite{shahapure2020cluster} guided our selection of the number of clusters $m$. 
We found that $m=5$ is a reasonable choice. This clustering approach for determining the number and characteristics of launch missions is shown in Figure \ref{fig:mission_cluster}. Note that historical launches are near circular, so we approximate the orbital radius as the semi-major axis.


The weight assigned to each mission type is determined by the proportion of data points in each cluster relative to the total data set. With the predicted total demand and derived weights, we estimate the demand $k_j$ for each mission type $j$. In terms of the capacity of a spaceport, we define $P$ as the maximum capacity of each spaceport and, in this paper, we assume a spaceport can launch at most one rocket a week. With the total launch demand and the capacity of a spaceport, the number of spaceports required is derived and denoted as $K$. Therefore, the demand constraints of the SPFLP can be formulated as 
\begin{subequations}
\begin{alignat}{2}
\sum_{i=1}^N x_{i} & = K \label{eq:demand_1a}, \\
\sum_{i=1}^N y_{i,j} & \geq k_{j}, \qquad\forall j = 1, \dots ,m\label{eq:demand_1b}, \\
\sum_{j=1}^m y_{i,j} & \leq P, \;\qquad\forall i = 1, \dots ,N \label{eq:demand_1c},
\end{alignat}
\label{eq:demand_const}
\end{subequations}
\begin{table}
\begin{center}
\begin{tabular}{|c|c c c c|}
\hline
Mission type & $a_\ell$ (km) & $\mathbbm{i}_{\ell}$ ($\degree$)& Weight (\%)& Demand \\
\hline
Type 1  &  6694.72 & 54.34  & 42\% & 116\\
\hline
Type 2  &  6990.01 & 43.30 & 6\% & 17\\
\hline
Type 3  &  7187.32 & 97.78 & 18\%& 48\\
\hline
Type 4  &  6893.20 & 95.91 & 31\% & 84 \\
\hline
Type 5  &  7649.00 & 75.46 & 3\%& 8\\
\hline
\end{tabular}

\label{tab:mission_clustering}
\end{center}
\caption{Mission type characteristics after identification through clustering.}
\end{table}

\noindent
where $N$ is the total number of counties and $x_i$ is a binary decision variable indicating whether or not to build a spaceport at county $i$. Constraints \eqref{eq:demand_1a} and \eqref{eq:demand_1b} ensure that the total launch demand in 2030 and the demands for each mission type are satisfied, respectively. Constraint \eqref{eq:demand_1c} guarantees that the total number of assigned mission for a spaceport will not exceed its maximum capacity. Note that we could have trivially enforced different per-spaceport capacities (i.e., different $P_i$'s). 

As mentioned previously, we introduce coupling constraints between spaceport placement decisions ($x_i$ decision variables) and mission allocation decisions ($y_{i, j}$ decision variables). 
Specifically, when a candidate county is \emph{not} selected, then there should be \emph{no} missions assigned to it. Conversely, when a location is selected for placement of a spaceport, at least one mission of a certain mission type must be assigned to it. Explicitly, these coupling constraints can be written as
\begin{subequations}
\begin{alignat}{2}
\frac{1}{KP} \sum_{j=1}^m y_{i,j} &\leq x_i, \qquad\qquad\forall i = 1, \dots ,N \label{eq:coup_1a},\\
x_i &\leq \sum_{j=1}^m y_{i,j}, \qquad\forall i = 1, \dots ,N \label{eq:coup_1b}.
\end{alignat}
\label{eq:coup}
\end{subequations}


\subsubsection{Spatial dispersal of spaceports}
We wish to guard against the scenario where the optimal solution selects a set of counties in close geographic proximity to each other. This aligns with, e.g., US Space Force objectives to build resilient space architectures \cite{resilient_space}, since spaceports in close proximity may be impacted by the same disruptive event (e.g., hurricanes, earthquakes, etc.) To ensure continuity of launch capabilities in the event of a natural disaster, each pair of selected spaceports must be adequately separated. We denote this minimum separation distance by $D$. Note that a readily implementable extension is to tailor $D$ according to the geographic extent of a given natural disaster. For simplicity (and because it does not affect the model formulation), we assume a constant minimum separation distance.

We denote by $d_{i, i'}$ the geographic distance between two candidate counties $i$ and $i'$. 
To enforce spatial dispersion, 
we need to ensure $d_{i, i'} \geq D$ for all pairs of selected candidate counties $i$ and $i'$. 
To model this logic in the SPFLP, we introduce a binary variable $z_{i, i'}$ for each pair of candidate counties and apply the Big $M$ method. This gives the following set of constraints
\begin{subequations}
\begin{align}
D - M\left(1 - z_{i, i'}\right) \leq d_{i, i'}, \forall i, i' = 1, \dots ,N; \; i \neq i' \label{eq:dist_1a}, \\
d_{i, i'} - Mz_{i, i'} \leq D, \forall i, i' = 1, \dots ,N; \; i \neq i' \label{eq:dist_1b},\\
x_{i} + x_{i'} \leq 1 + z_{i, i'}, \forall i, i' = 1, \dots ,N; \; i \neq i' \label{eq:dist_1c},
\end{align}
\label{eq:distance}
\end{subequations}
\noindent
where $M$ is a very large, positive number. 

\subsubsection{Complete SPFLP formulation}

The complete formulation of the SPFLP model is
\begin{subequations}
\begin{alignat}{2}
\min_{x, y, z} \quad &\sum_{i = 1}^N \left(C_{T_i} + C_{O_i} \right)x_{i} + \sum_{i=1}^N \sum_{j=1}^m \left(C_{L_{i,j}} + C_{R_{i,j}} \right)y_{i,j} \label{eq:basic_obj},\\
\textrm{s. t.} \quad & \text{Constraints } \eqref{eq:demand_1a} - \eqref{eq:dist_1c}, \\
&x_{i}, z_{i, i'} \in \{0,1\}, \forall i, i' = 1, \dots ,N; \; i \neq i' \label{eq:basic_j},\\
&y_{i, j} \in \mathbb{N}_{\geq0}, \forall i = 1, \dots, N; \; \forall j = 1, \dots, m \label{eq:basic_k}.
\end{alignat}
\label{eq:basic_flp}
\end{subequations}
\noindent
$C_{T_i}$ and $C_{O_i}$ denote the  transportation and operation cost associated with a spaceport at each candidate county $i$. $C_{L_{i,j}}$ represents the cost of launching mission type $j$ from a spaceport in candidate county $i$  (Section \ref{sssec:launc_traj_costs}), while $C_{R_{i,j}}$ refers to the aircraft rerouting cost  associated with launching mission type $j$ from a spaceport in candidate county $i$. 
\section{Experiment and sensitivity analysis}  \label{sec:experiment}

In this section, we present the selected counties and mission allocations, alongside an analysis of their impact on NAS flight operations. Additionally, a sensitivity analysis is conducted to assess impacts to spaceport selection and air traffic under different SPFLP configurations. In this section, the terms \emph{configuration} and \emph{scenario} are used as follows:
\begin{itemize}
    \item \emph{Scenario}: A cost weighting scenario is a distribution of weights placed on operation cost, transportation cost, mission launch cost, and aircraft rerouting cost. The four scenarios tested are explained in Section  \ref{subsec:experiment_setup}.
    \item \emph{Configuration}: A unique set of parameters used to run the SPFLP model. For each set of parameters, we choose one hazard zone size, one cost weighting scenario, and one air traffic volume level.
\end{itemize}





\subsection{Experiment setup}   \label{subsec:experiment_setup}
To evaluate the impact of spaceports on NAS flight operations, we focus on estimating flight rerouting costs under different spaceport location results from the SPFLP. We calculate flight rerouting costs by multiplying the number of flights within the hazard zones of each chosen spaceport with the unit cost of rerouting a flight (\$293 per flight as estimated in \cite{JohnHansman2018}). In our analysis, we used OpenSky flight trajectory data \cite{schafer2014bringing} from two days with differing traffic levels (normal air traffic level day on Feb. 12, 2022 with approximately 745,000 flights and a heavy air traffic level day on Nov. 24, 2022 with more than 940,000 flights) to represent distinct air traffic levels.

Moreover, we ran the model under several different configurations to evaluate the sensitivity of results when parameters are changed. The adjusted parameters are air traffic levels (high and low), hazard zone spreads ($\xi = 5\degree, 7.5\degree, 10\degree$), and cost weighting scenarios. 
The four cost weighting scenarios S1, S2, S3, and S4 are:
\begin{itemize}
    \item[S1:] Transportation cost, operation  cost, mission launch cost, and flight rerouting cost are all given equal weight. 
    \item[S2:] Launch cost is weighted $10\times$ more than other costs. 
    \item[S3:] Rerouting cost is weighted $10\times$ more than other costs. 
    \item[S4:] Transportation cost and operation cost are weighted $10\times$ more than other costs.
\end{itemize}



In every configuration, the SPFLP selects six optimal spaceports which satisfy all constraints while minimizing total costs. 
We set $D = 300$ miles for all model runs; 
Our selection of $300$ miles was guided by the fact that this is the average width of a hurricane \cite{hurricane_size}. 
An example of selected spaceports under hazard zone size $\xi = 5\degree$, low air traffic, but with different cost weighting scenarios is shown geographically in Figure \ref{fig:d10_traj1}. 

\subsection{Sensitivity analysis and discussion}
\subsubsection{Sensitivity to sizes of hazard areas}
\begin{figure*}[!htbp]
    \centering
    \setlength\fboxsep{0pt}
    \setlength\fboxrule{0pt}
    \fbox{\includegraphics[width=7in]{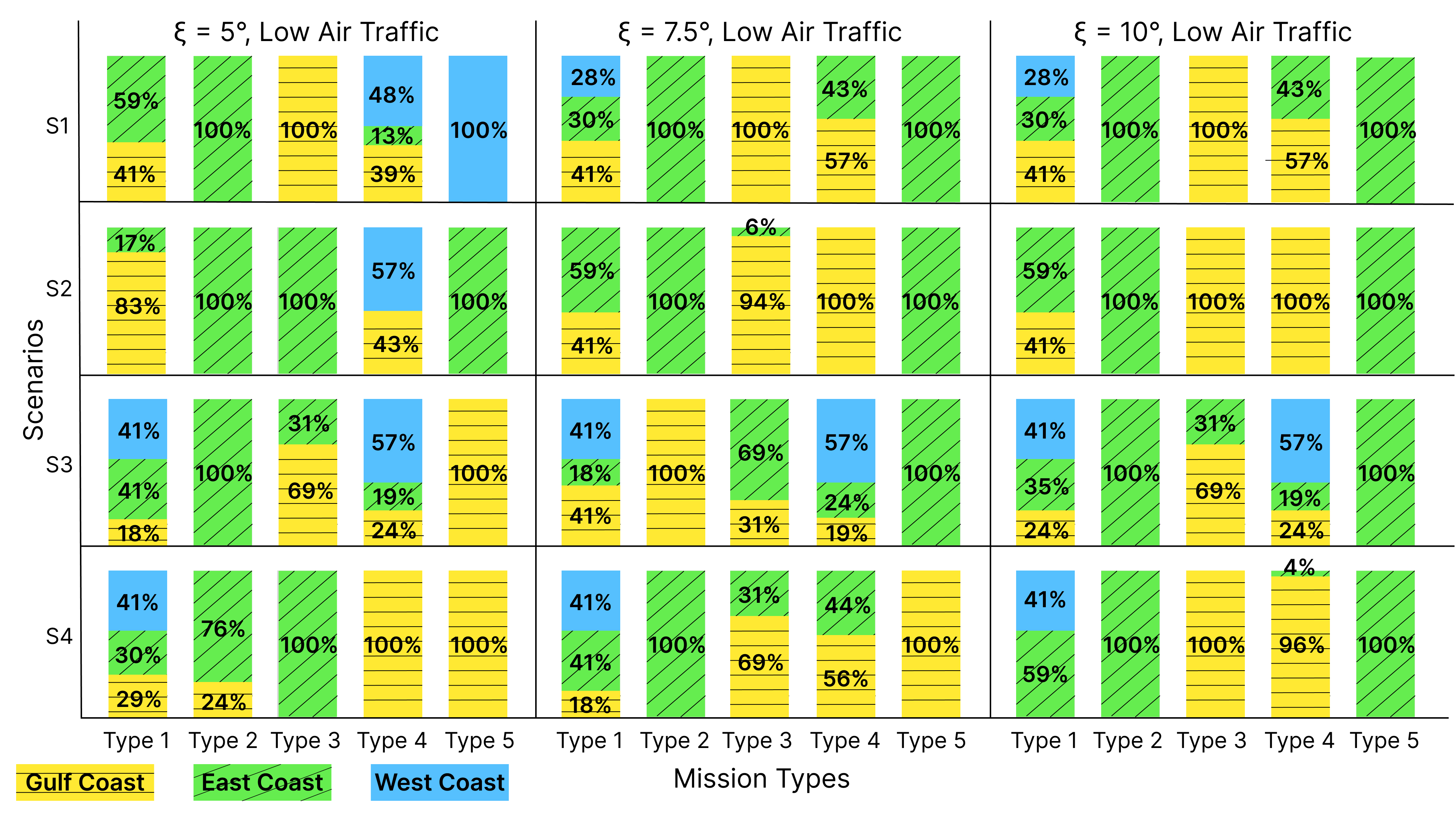}}
    \caption{Spaceport location planning results showing breakdown of launch mission type allocations among West Coast, East Coast, and Gulf Coast spaceport locations for low levels of air traffic.}
    \label{fig:low_traffic_region_selections}
\end{figure*}

The SPFLP analysis under low air traffic and Scenario 1 cost weighting, with varying hazard zone sizes depicted in Figure \ref{fig:low_traffic_region_selections}, reveals that West Coast spaceport mission allocations are significantly influenced by hazard zone sizes. At $\xi = 5\degree$, West Coast spaceports received 48\% of Type 4 and 100\% of Type 5 missions. For $\xi = 7.5\degree$ and $\xi = 10\degree$, these spaceports were allocated 28\% of Type 1 missions. Notably, Humboldt, CA, was selected for $\xi = 5\degree$, but the choice shifted to Del Norte, CA, for larger $\xi$ values, as seen in Figure \ref{fig:d10_traj1}. This shift is likely due to increased hazard sizes making Humboldt's launch trajectories non-compliant with population avoidance criteria, diminishing its suitability for launches. Thus, the West Coast's adaptability to population constraints and airline rerouting costs escalates with hazard zone size.

\subsubsection{Sensitivity to the mission launch cost}
\begin{figure}[!htbp]
    \centering
    \subfloat[\label{subfig:d10_traj1_1111}$\text{Scenario 1}$]{\includegraphics[width=4.4cm]{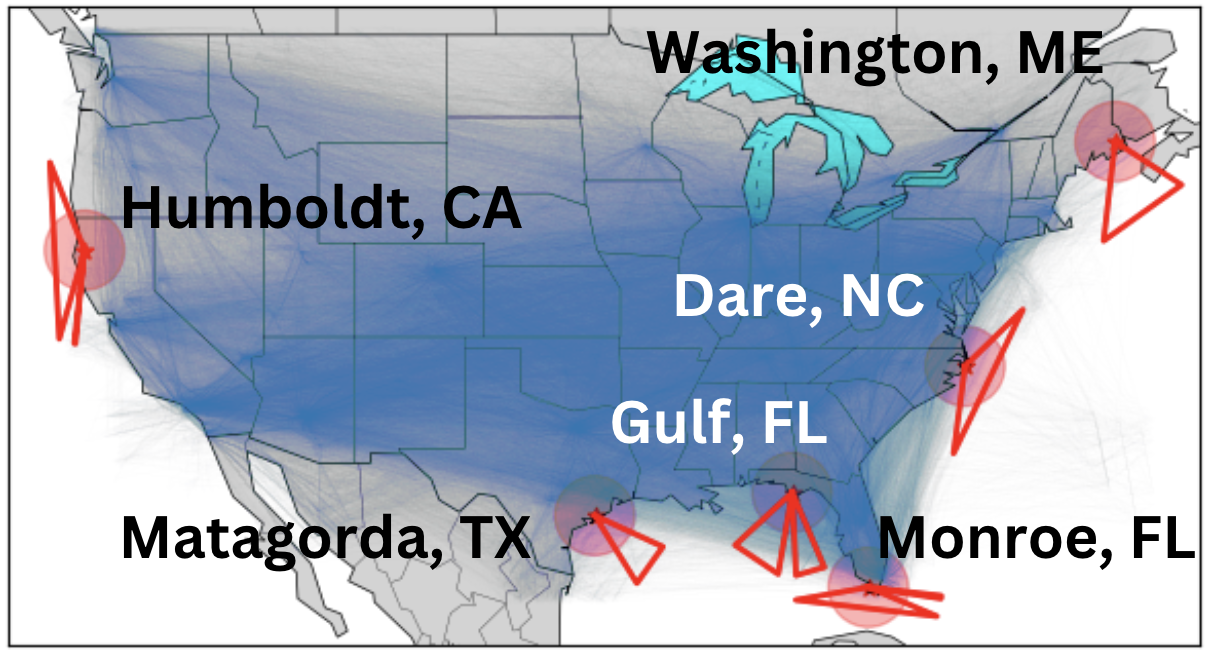}}\hfill
    \subfloat[\label{subfig:d10_traj1_11110} $\text{Scenario 2}$]{\includegraphics[width=4.4cm]{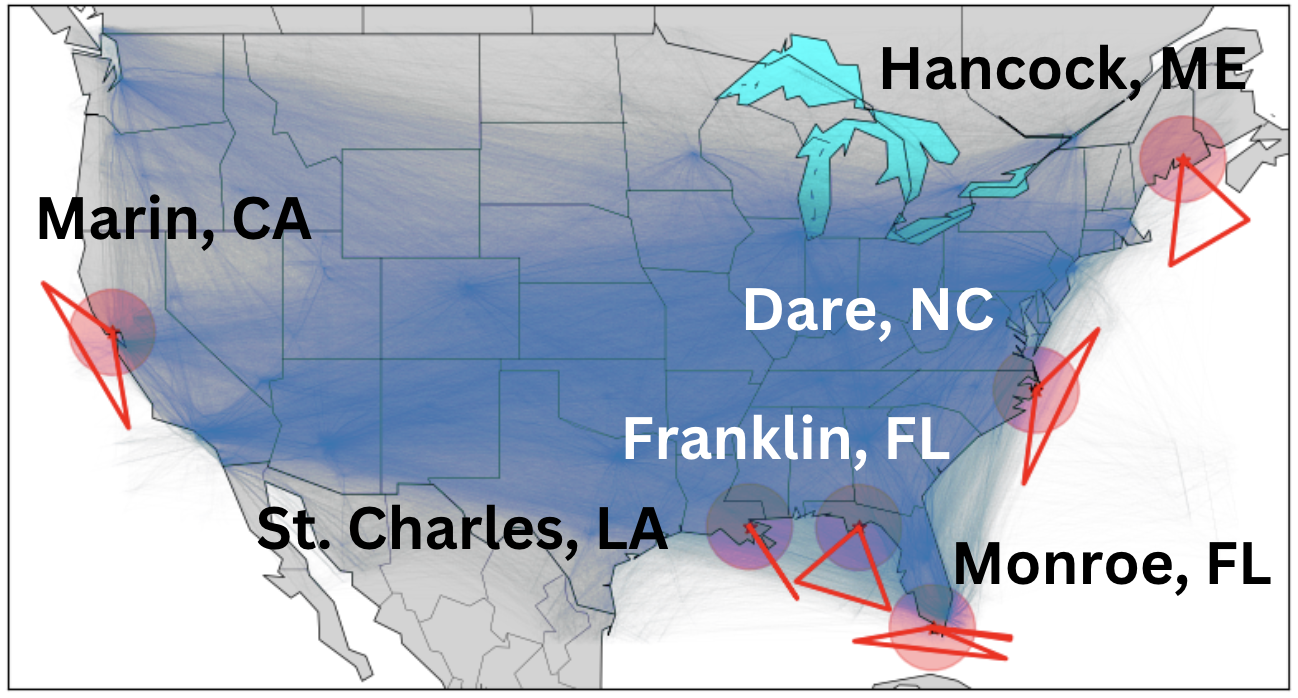}}\vfill
    \subfloat[\label{subfig:d10_traj1_11101}$\text{Scenario 3}$]{\includegraphics[width=4.4cm]{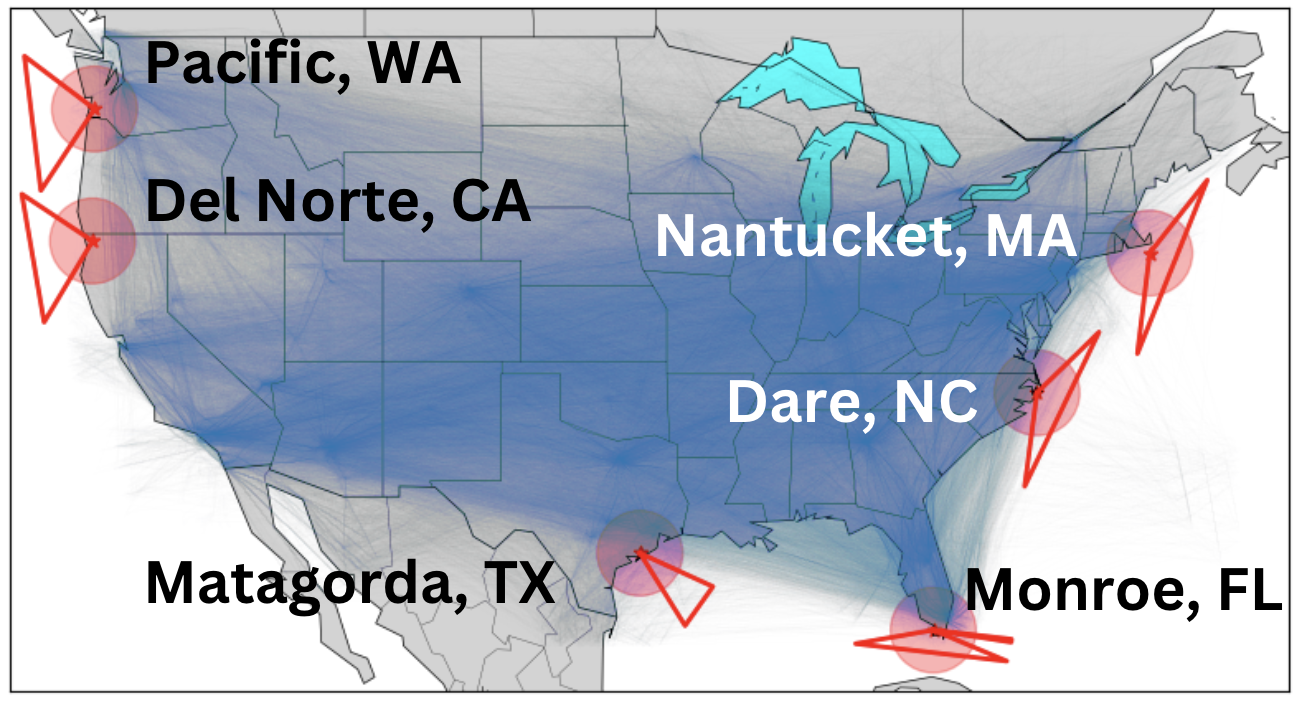}}\hfill
    \subfloat[\label{subfig:d10_traj1_101011}$\text{Scenario 4}$]{\includegraphics[width=4.4cm]{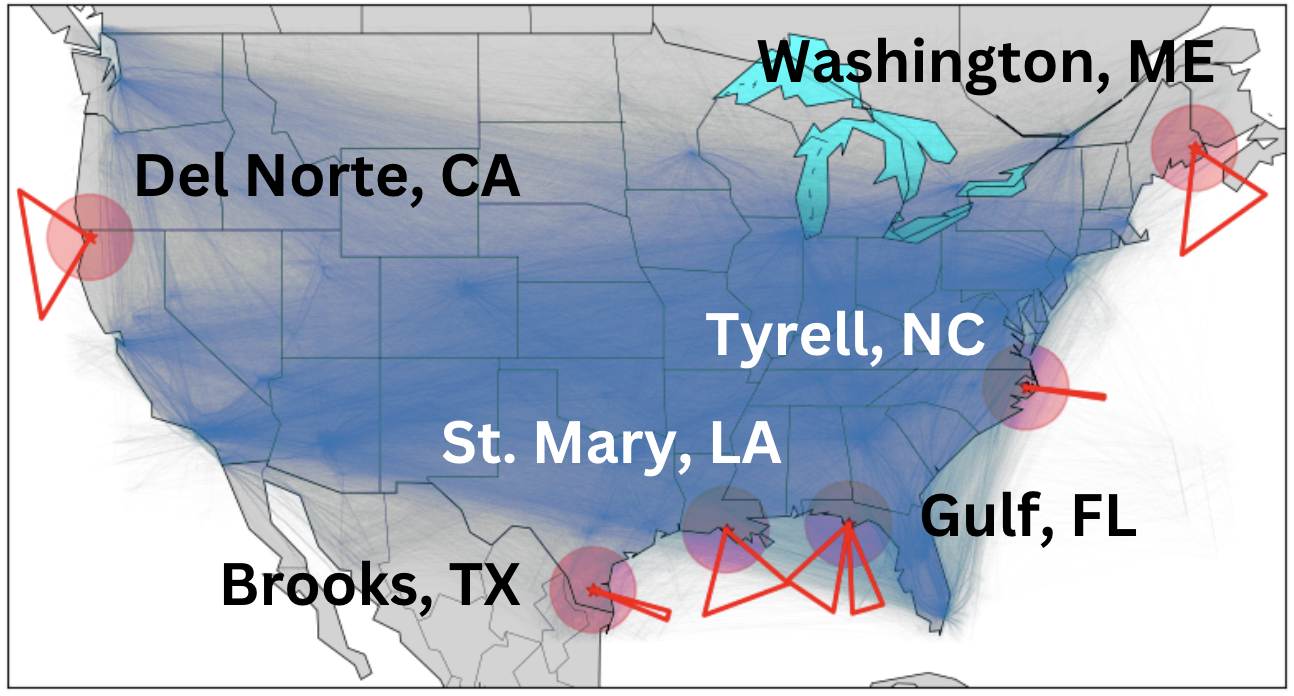}}\vfill
    \caption{$\xi = 5^{\circ}$ with low air traffic. Red triangles indicate the hazard zones, and each selected spaceport location has a red disc with radius $D = 300$ miles.}\label{fig:d10_traj1}
\end{figure}

In Scenario 2, there is a notable increase in Type 1 missions launched from Gulf Coast spaceports and a decrease for West Coast spaceports. As depicted in Figure \ref{subfig:d10_traj1_11110}, launch angles for the Gulf Coast are approximately $180\degree$, while those for the West Coast are closer to $270\degree$. Utilizing the launch cost methodology described in Section \ref{sec:methodology}, we find that the launch cost for Type 1 missions is highest at a launch angle of $270\degree$, and lowest for launch angles around $150\degree$. Consequently, Gulf Coast launches dominate in Scenario 2 when launch costs are dominant. Scenario 3 illustrates a shift towards the West Coast due to dominant rerouting costs and denser air traffic trajectories around the Gulf Coast. Therefore, Gulf Coast spaceports are preferable for minimizing launch costs as compared to their West Coast counterparts
\subsubsection{Sensitivity to transportation and operation costs}
In Scenario 4, the emphasis on transportation and operation costs results in missions being served by spaceports on all coasts. When transportation and operation costs have greater significance, there is no clear bias for specific coasts: However, counties characterized by lower median house values and reduced average traffic times are favored. 
\subsubsection{Sensitivity to the air traffic level} \label{sssec:sensitivity_airtraffic}
High and low air traffic analyses yield comparable mission distributions, with Types 2 and 5 missions preferring the East Coast, and Type 3 the Gulf Coast, across most scenarios. Type 1 and 4 missions show varied allocations based on scenario and hazard zone sizes, but maintain consistent trends across traffic levels. Table \ref{tab:rerouting_costs} indicates that rerouting costs generally rise with increased air traffic. In contrast,  Scenario 3 presents lower rerouting costs at higher traffic levels due higher rerouting cost weight.
\subsubsection{Mission type preferences}
From our results and sensitivity analyses, we observe that the East Coast is favored for launching Types 2 and 5 missions, the Gulf Coast for Type 3 and 4 missions, and the West Coast is less suitable for Type 1 and 2 missions due to higher launch costs.

%
%

\begin{figure*}[!htbp]
    \centering
    \setlength\fboxsep{0pt}
    \setlength\fboxrule{0pt}
    \fbox{\includegraphics[width=7in]{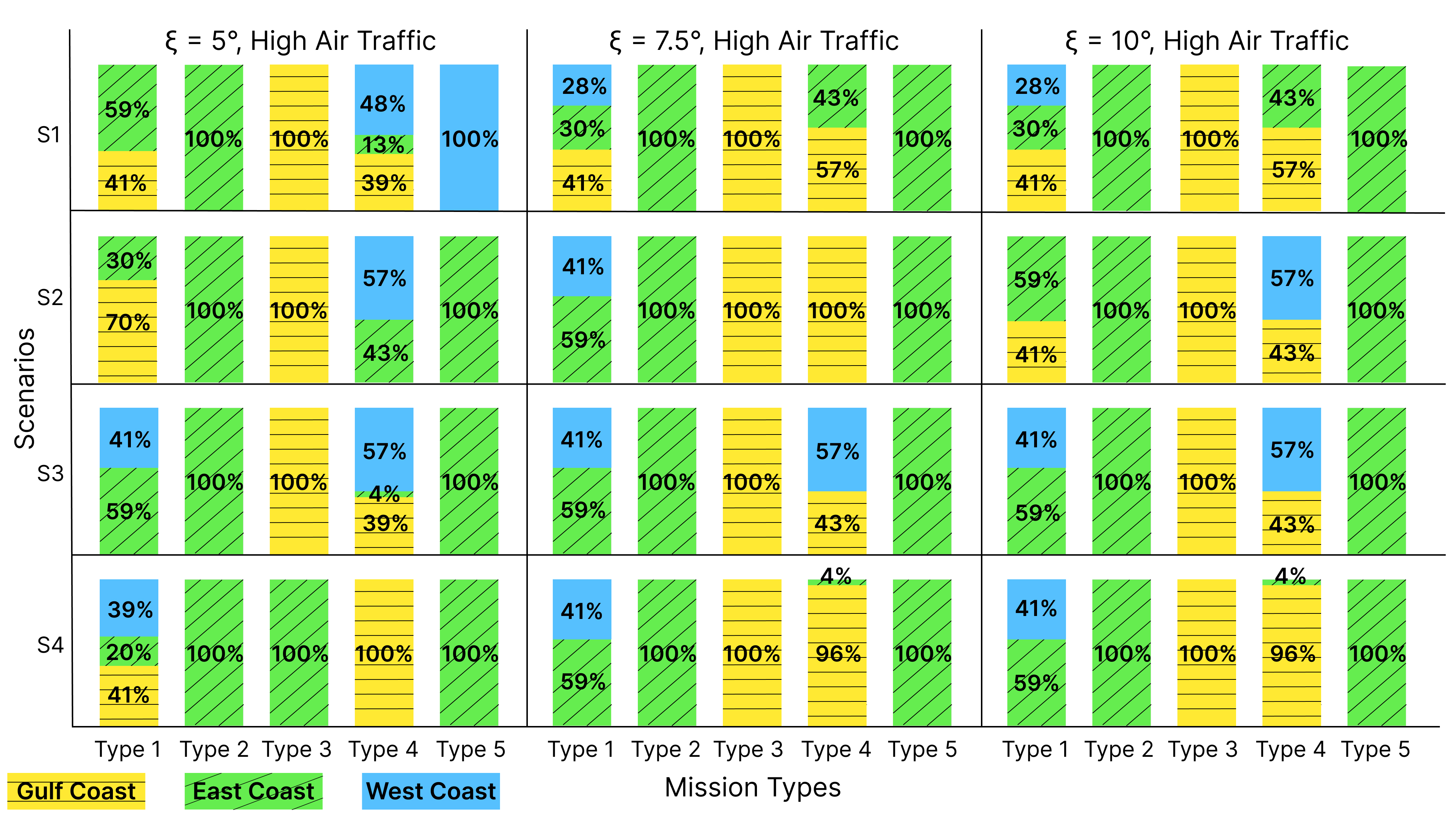}}
    \caption{Spaceport location planning results showing breakdown of launch mission type allocations among West Coast, East Coast, and Gulf Coast spaceport locations for high levels of air traffic.}
    \label{fig:high_traffic_region_selections}
\end{figure*}

\begin{table*}
\centering
\caption{Spaceport location planning results: Costs in millions of US dollars (\$) for different SPFLP configurations. Note that Scenario 3 ($*$) results are counterintuitive, as discussed in Section \ref{sssec:sensitivity_airtraffic}.}
\begin{center}
\begin{tabular}{|c|c c|c c|c c|}
\hline
\multirow{2}{*}{Scenario} & \multicolumn{2}{c|}{$\xi = 5^{\circ}$} & \multicolumn{2}{c|}{$\xi = 7.5^{\circ}$}  & \multicolumn{2}{c|}{$\xi = 10^{\circ}$}\\
& Low traffic & High traffic & Low traffic & High traffic & Low traffic & High traffic \\
\hline
$\textbf{S}_1$ \textcolor{white}{($*$)}  &  18.6 & 18.8 & 19.0 & 18.8 & 20.3 & 20.5 \\
\hline
$\textbf{S}_2$ \textcolor{white}{($*$)}  &  27.5 & 28.9 & 24.4 & 26.9 & 25.4 & 28.0 \\
\hline
$\textbf{S}_3$ ($*$)   & 9.8 & 8.4 & 10.1 & 8.8& 10.6 & 9.5 \\
\hline
$\textbf{S}_4$ \textcolor{white}{($*$)}  &  19.6 & 20.3 & 21.2 & 22.5 & 20.6 & 21.9 \\
\hline
\end{tabular}
\end{center}
\label{tab:rerouting_costs}
\end{table*}
\section{Concluding remarks}    \label{sec:conclusion}
We introduce a framework designed to rigorously select potential spaceport locations and devise mission assignment plans that meet future space launch demand forecasts. Our spaceport facility location planning model considers a wide range of costs and objectives to select potential spaceport locations. Our study contains several limitations and strong assumptions: Some were necessitated by the lack of publicly available data from space launch operators. Other limitations and assumptions can be addressed or relaxed in future work: For example, one could consider more precise cost assessments, particularly for transportation and operation costs. Additionally, more sophisticated robust optimization-based models could be used to plan resilient spaceport networks in the face of natural disasters and other potential disruptions.

\bibliographystyle{IEEEtran} 
\small{
\bibliography{main.bib}
}


\end{document}